\documentclass[conference]{IEEEtran}
\IEEEoverridecommandlockouts

\usepackage[utf8]{inputenc}
\usepackage{cite}
\usepackage{amsmath,amssymb,amsfonts,amsthm}
\usepackage{graphicx}
\usepackage{textcomp}
\usepackage{xcolor}
\usepackage{bm}
\usepackage{array}
\usepackage{booktabs}
\usepackage{url}

\newtheorem{remark}{Remark}

\newcommand{\calV}{\mathcal{V}}
\newcommand{\umax}{u_{\max}}
\newcommand{\sat}{\operatorname{sat}}

\begin{document}

\title{Beyond Phase Reduction: Amplitude Collapse in Optimal Control of Coupled Oscillators}

\author{%
\IEEEauthorblockN{Faranak Rajabi}
\IEEEauthorblockA{\textit{Dept.\ of Mechanical Engineering}\\
\textit{University of California, Santa Barbara}\\
faranakrajabi@ucsb.edu}
\and
\IEEEauthorblockN{Fr\'ed\'eric Gibou}
\IEEEauthorblockA{\textit{Dept.\ of Mechanical Engineering}\\
\textit{University of California, Santa Barbara}\\
fgibou@ucsb.edu}
\and
\IEEEauthorblockN{Jeff Moehlis}
\IEEEauthorblockA{\textit{Dept.\ of Mechanical Engineering}\\
\textit{University of California, Santa Barbara}\\
moehlis@ucsb.edu}
}

\maketitle
\thispagestyle{plain}
\pagestyle{plain}

\begin{abstract}
We solve the four-dimensional Hamilton-Jacobi-Bellman (HJB)
equation for two diffusively coupled Stuart-Landau-like
oscillators to obtain full-state optimal feedback control.
A sweep over coupling strength reveals a
sharp change in the numerically optimal strategy: below a threshold coupling
value, the controller steers the phase difference toward
anti-phase while keeping both oscillators near the limit cycle,
as reduced-order models would suggest. Above this threshold, the HJB solution changes qualitatively; the controller
transiently collapses one oscillator's amplitude to near zero,
thereby enabling large phase repositioning near the origin before rebuilding its amplitude. 
Direct gradient-based and stochastic optimization do not recover this lower-cost collapse trajectory from the initializations considered, suggesting that it occupies a region of the control landscape that is difficult to access by direct search.
A joint sweep over
nonisochronicity and coupling shows that collapse can occur even for an isochronous oscillator: phase repositioning near the origin can favor an off-cycle strategy.  Nonisochronicity provides an additional energetic benefit through a phase-velocity surplus at small amplitude, quantitatively accounting for the observed reduction in control cost.
Comparisons with uncoupled and coupled phase-reduced baselines show that phase models 
become increasingly inaccurate and cost significantly more energy for strong coupling.
Results for coupled Van~der~Pol oscillators further demonstrate that exploitation of off-cycle dynamics
is not specific to the Stuart-Landau-like oscillators.
\end{abstract}


\section{Introduction}\label{sec:intro}

Coupled limit-cycle oscillators arise throughout science and
engineering, and controlling their phase relationships has direct
clinical and infrastructure consequences. In deep brain
stimulation for Parkinson's disease, pathological synchrony of
neural populations drives motor symptoms; current open-loop
protocols (${\gtrsim}130$\,Hz pulse trains) are effective but
energy-intensive, requiring surgical battery replacement every
3 to 5 years~\cite{tass1999,brocker2017}. Energy-optimal control
could extend device lifetimes and enable emerging closed-loop
paradigms. In power systems, grid-forming inverters increasingly
employ virtual oscillator control based on coupled
Hopf-normal-form dynamics~\cite{johnson2016,seo2019}, where
desynchronization can cascade to system-wide
failure~\cite{dorfler2012}. In both domains, the fundamental
problem is identical: steer coupled oscillators toward a desired
phase relationship while minimizing energy.

Phase reduction~\cite{kuramoto1984,ermentrout2010} is the
standard tool for oscillator control, projecting each oscillator
onto a single phase variable evolving on the limit cycle. This
has enabled elegant optimal control
solutions~\cite{moehlis2006,zlotnik2013,monga2019}, but
it constrains the dynamics to a neighborhood of the cycle by construction. When
perturbations or coupling become stronger, pure phase descriptions
break down and higher-order or phase-amplitude descriptions
become
necessary~\cite{wilson2019,wilson2020b,franci2012,nicks2023}.
Recent work has pushed reduced descriptions into stronger-input
regimes using adaptive phase-amplitude coordinates and dynamic
programming on the reduced
model~\cite{wilson2020adapt,wilson2021b,dewanjee2024}.

Here we investigate whether full-state optimal control can select qualitatively different mechanisms once coupling is strong enough that remaining near the limit cycle is costly. Solving the four-dimensional HJB equation reveals two regimes. At weak coupling, the optimal trajectories remain near the limit cycle and closely resemble phase-based control. At stronger coupling, the computed HJB solution instead drives one oscillator to very small amplitude before returning it to the limit cycle. Near the origin, the oscillator's phase becomes highly sensitive to state perturbations, enabling large phase repositioning and, for nonisochronous oscillators, can additionally exploit a higher angular velocity at small amplitude. Thus, the full-state calculation identifies an off-cycle control mechanism that is absent from standard phase-reduced descriptions and becomes favorable over a well-defined range of parameters.


\section{Related Work}\label{sec:related}

\textbf{Phase-based optimal control.}
Classical phase reduction~\cite{kuramoto1984} projects each
oscillator to a single angular variable, enabling elegant
minimum-energy and minimum-time control
formulations~\cite{moehlis2006,zlotnik2013}.
Phase models can incorporate coupling through interaction
functions~\cite{wilson2022}, and recent work extends
optimal control to stronger inputs using adaptive
phase-amplitude coordinates~\cite{wilson2020adapt,wilson2021b,monga2020}
and phase-based control of strongly perturbed
oscillators~\cite{dewanjee2024,namura2024,takata2021}.
By construction, these methods constrain the state to lie on or
near the limit cycle. 
Several papers document that standard phase reduction can miss
important dynamics under stronger coupling or forcing, including
bifurcations and amplitude-mediated
effects~\cite{franci2012,wilson2019,wilson2020b,nicks2023}.
Isostable coordinates~\cite{wilson2016} provide a principled way
to capture transverse dynamics, and adaptive
frameworks~\cite{wilson2020adapt,dewanjee2024} explicitly
aim to relax the small-perturbation limitation. 
Although these techniques have not been used to capture excursions close to the origin, they
suggest that a full-state computation may reveal qualitatively
new phenomena. Direct trajectory optimization avoids reduction altogether but
returns local optima that depend on initialization; we compare
against gradient-based and evolutionary formulations in
Sec.~\ref{sec:direct}.

\textbf{HJB for oscillator control.}
HJB methods have been applied to oscillator control problems in
neuroscience, including driving oscillators to phaseless sets for
phase randomization~\cite{danzl2009}, seizure
termination~\cite{wilson2014seizure}, and minimum-energy
desynchronization~\cite{nabi2013,wilson2014}. The work
of~\cite{danzl2009} is the closest precedent, as it also uses
HJB to exploit off-cycle geometry. However, that paper targets
phaseless-set randomization for a single oscillator rather than
controlled phase repositioning in a coupled pair. More broadly,
control and synchronization of coupled oscillator networks have
been studied with objectives beyond pure phase
models~\cite{skardal2016,skardal2022,salfenmoser2024,trummel2023}.
Related multi-agent work addresses predictive consensus and adaptive
formation control~\cite{rezaei2023,norouzi2025}. To the best of our
knowledge, none of this work reports a globally optimal off-cycle
amplitude-collapse mechanism for diffusively coupled limit-cycle
oscillators.

HJB methods return
a global optimum over all
admissible controls, with existence and uniqueness of the value
function supplied by viscosity solution theory~\cite{crandall1983},
whereas direct methods return a local optimum with no certificate
that a lower-cost basin was missed. This matters beyond the dimension
solved: certificates computed in a reduced representation transfer to
the full system only under conditions on the geometry the reduction
preserves~\cite{lutkus2025,nakamura2025}, and such transfer
presupposes a trustworthy reduced solution. The present results show
that a reduced description can be wrong without being visibly wrong,
since phase reduction returns a plausible on-cycle strategy while a
cheaper off-cycle strategy exists outside its domain of validity.

\section{Problem Formulation}\label{sec:problem}

\subsection{Two Coupled Stuart-Landau-Like Oscillators}

We consider two diffusively coupled oscillators whose dynamics
take the Stuart-Landau-like form:
\begin{align}
\dot{x}_i
  &= a x_i - \omega y_i - a x_i^3 - (a + 2b)x_i y_i^2
     + c(x_j - x_i) + u_i,
     \label{eq:SLx}
\\
\dot{y}_i
  &= \omega x_i + a y_i - a y_i^3 + (2b - a)x_i^2 y_i,
     \label{eq:SLy}
\end{align}
for $i \in \{1,2\}$, $j \neq i$, with $a > 0$, $\omega > 0$,
$b \in [0,\omega)$, $c \ge 0$, and
$u_i \in [-\umax,\umax]$. Control and coupling enter only the
$x$-channel, modeling electrode-based actuation and gap-junction
coupling in neural applications~\cite{nabi2013}.

\begin{figure}[t]
\centering
\includegraphics[width=\columnwidth]{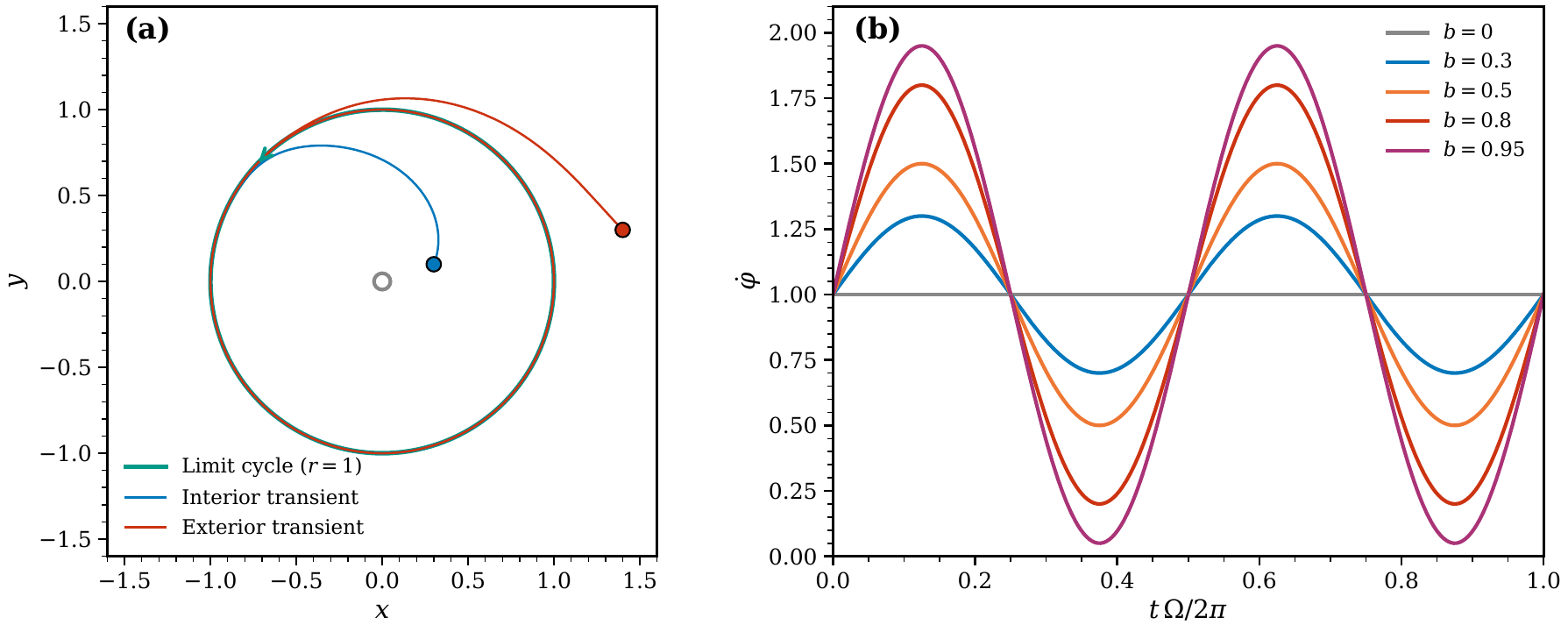}
\caption{Uncoupled Stuart-Landau-like dynamics.
(a)~Transient trajectories from interior and exterior initial
conditions spiral onto the unit-circle limit cycle.
(b)~Instantaneous angular velocity
$\dot\phi = \omega + br^2\sin 2\phi$ on the limit cycle for
several nonisochronicity values. At $b = 0$ the oscillator is
isochronous ($\dot\phi \equiv 1$). As $b$ increases,
$\dot\phi$ oscillates more strongly and the time-averaged
frequency $\Omega = \sqrt{\omega^2 - b^2}$ drops below the
origin frequency $\omega$. This gap is the phase velocity
surplus exploited by the collapse strategy
(Sec.~\ref{sec:mechanism}).}
\label{fig:autonomous}
\end{figure}

We refer to these as ``Stuart-Landau-like'' because of the extra
$\sin 2\phi$ term that appears in the polar dynamics, which does not appear in the standard Stuart-Landau normal form. In polar
coordinates $(x_i,y_i) = r_i(\cos\phi_i,\sin\phi_i)$, the
uncoupled dynamics are
\begin{equation}
\dot r = ar(1-r^2), \qquad
\dot\phi = \omega + br^2\sin 2\phi.
\label{eq:polar}
\end{equation}
Fig.~\ref{fig:autonomous} illustrates these dynamics for the
uncoupled system. 
In addition to a stable in-phase state, there is an anti-phase solution with
$x_2(t) = -x_1(t)$ and $y_2(t) = -y_1(t)$, providing a natural
target state for the control problem. The anti-phase solution is stable for very weak coupling $c \lesssim 0.015$.  For larger coupling strength, including those considered below, the transverse Floquet multiplier for the anti-phase solution exceeds unity, and trajectories starting near this solution 
return to the in-phase solution. The objective is therefore to reach an unstable
target orbit within a finite horizon, which is why the terminal cost
must anchor the radii explicitly. The unit-circle limit
cycle has period $T = 2\pi/\Omega$ where
$\Omega = \sqrt{\omega^2 - b^2}$ and Floquet exponent
$\lambda = -2a$. We fix $a = \omega = 1$, $b = 0.8$,
$\umax = 1$ unless stated otherwise. Setting $a = \omega = 1$ is a choice of units, not a
restriction: rescaling time fixes $\omega = 1$ and rescaling
amplitude places the cycle at $r = 1$, fixing $a = 1$.
Nonisochronicity is bounded by $b < \omega$ so that
$\Omega = \sqrt{\omega^2-b^2}$ is real, and the period diverges as
$b \to 1$, which is why the sweeps stop at $b = 0.9$. We will refer to $b$ as the
nonisochronicity parameter.

\begin{remark}[Symmetry]\label{rem:symmetry}
The system is $\mathbb{Z}_2$-equivariant under oscillator swap~\cite{gss}.
The HJB equation inherits this symmetry; we verify it
\emph{a~posteriori} to within discretization error. The initial
condition $z_0  \equiv (x_1(0),y_1(0),x_2(0),y_2(0)) =  (1,0,\cos 0.1,\sin 0.1)$ breaks this symmetry;
the controller selects which oscillator to collapse based on the
state.
\end{remark}

\subsection{Optimal Control Objective}

The cost functional is
\begin{equation}
J(z_0) = \int_0^{T_f} R(u_1^2 + u_2^2)\,dt
         + \frac{\gamma}{2}\,d^2(z(T_f)),
\label{eq:cost}
\end{equation}
with $R = 10$, $\gamma = 1000$, $T_f = 2T$, and terminal cost
\begin{equation}
d^2(z) = (x_1{+}x_2)^2 + (y_1{+}y_2)^2
         + (r_1{-}1)^2 + (r_2{-}1)^2
\label{eq:orbit_target}
\end{equation}
targeting the entire anti-phase orbit. 
The cost (\ref{eq:cost}) is expected to depend on the initial condition $z_0$.
Since $d^2$ depends on
$(x_1{+}x_2, y_1{+}y_2)$ and radii $r_i$ but not the absolute
phase, it vanishes for any anti-phase pair
$(r\,e^{i\phi}, r\,e^{i(\phi+\pi)})$ with $r = 1$, regardless
of $\phi$. The $(r_i{-}1)^2$ terms prevent the trivial solution
$z = 0$ from achieving zero terminal cost.

\section{HJB Solution Method}\label{sec:hjb}

Let $\calV(z,\tau)$ denote the value function in backward time
$\tau = T_f - t$, where $z = (x_1,y_1,x_2,y_2) \in \mathbb{R}^4$
is the full state. It satisfies the Hamilton-Jacobi-Bellman equation
\begin{equation}
\frac{\partial \calV}{\partial \tau}
+ H\!\left(z,\nabla_z \calV\right) = 0,
\qquad
\calV(z,0) = \frac{\gamma}{2}\,d^2(z),
\label{eq:HJB}
\end{equation}
where $d(z)$ denotes the distance to the target set. The Hamiltonian is
\begin{equation}
H(z,p) = -f(z)\cdot p + H^*(p_{x_1}) + H^*(p_{x_2}),
\end{equation}
with $p = \nabla_z \calV$. The convex dual $H^*$, corresponding to the
running cost $R u^2$ under $|u| \le \umax$, is
\begin{equation}
H^*(p) =
\begin{cases}
\dfrac{p^2}{4R}, & |p| \le 2R\umax, \\[6pt]
-R\umax^2 + \umax |p|, & \text{otherwise},
\end{cases}
\end{equation}
yielding the feedback control
\begin{equation}
u_i^* = \sat\!\left(-\frac{p_{x_i}}{2R},\,\umax\right),
\end{equation}
where
$\sat(v,\umax) = \max \{ -\umax,\min\{v,\umax\}\}$.

The PDE is solved on $[-1.5,1.5]^4$ using a uniform grid with
$G=64$ points per dimension. We employ a fifth-order WENO scheme
with local Lax-Friedrichs flux splitting and SSP-RK3 time stepping.
The solver extends CASL-HJX~\cite{rajabi2026caslhjx} to 4D with GPU
acceleration~\cite{caslhjx4d}, enabling full-state solutions in
approximately five minutes on a single NVIDIA A100.

Parameters $R=10$ and $\gamma=1000$ are selected from a sweep over
$R \in \{1,10,25,50,100\}$ and $\gamma \in \{500,1000,2000\}$ to
balance control effort and value-function regularity. All optimal
trajectories remain well within the computational domain
($\|z(t)\|_\infty < 1.2$), and enlarging the domain to
$[-2.0,2.0]^4$ changes the optimal cost by less than $0.3\%$.
Grid refinement to $G=80$ and $G=96$ changes the optimal cost by
less than $1\%$ at the couplings reported below. Wall times per
parameter set on a single A100 are $300$\,s at $G = 64$, $819$\,s
at $G = 80$ and $2397$\,s at $G = 96$.

\section{Main Results}\label{sec:results}

\subsection{Two Qualitatively Distinct Strategies}

\begin{figure}[t]
\centering
\includegraphics[width=\columnwidth]{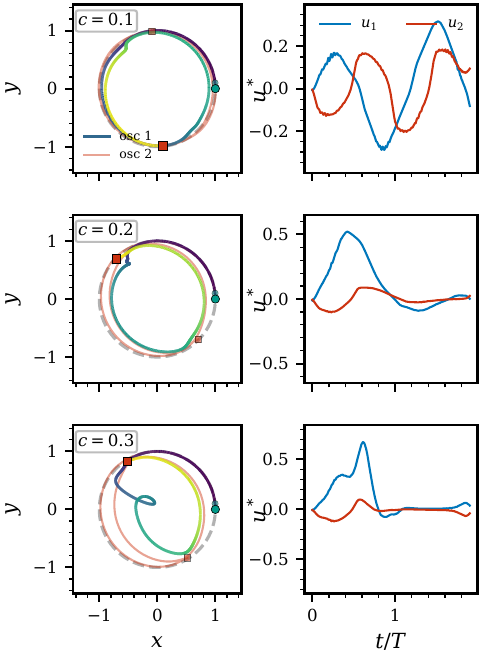}
\caption{Optimal trajectories for three coupling strengths.
Left column: phase portraits of oscillator~1 (color, graded by
time) and oscillator~2 (red) with the limit cycle shown dashed.
Green circle: initial state; red square: terminal state.
Right column: optimal feedback controls $u_1^*(t)$ and
$u_2^*(t)$. At $c = 0.10$, both oscillators remain near the
limit cycle and the controller applies small periodic
corrections. At $c = 0.20$, stronger coupling requires larger
control effort but the trajectory stays on-cycle. At $c = 0.30$,
the controller collapses oscillator~1 to $r_{1,\min} = 0.11$,
decoupling it from oscillator~2 and enabling free phase
repositioning before amplitude rebuilds naturally.}
\label{fig:dashboard}
\end{figure}

Fig.~\ref{fig:dashboard} shows HJB-optimal trajectories at
three coupling strengths, revealing two qualitatively distinct
control strategies. At weak coupling ($c = 0.10$), both
oscillators remain near the limit cycle throughout the control
horizon. The controller applies small perturbations to smoothly
steer the phase difference toward $180^\circ$. The trajectories
in the $(x_i, y_i)$ phase plane trace paths close to the unit
circle, and the minimum amplitude stays above
$r_{1,\min} > 0.85$. This is consistent with standard phase
reduction intuition: the oscillators stay near their limit cycle
while the controller adjusts the timing.

\begin{figure}[t]
\centering
\includegraphics[width=\columnwidth]{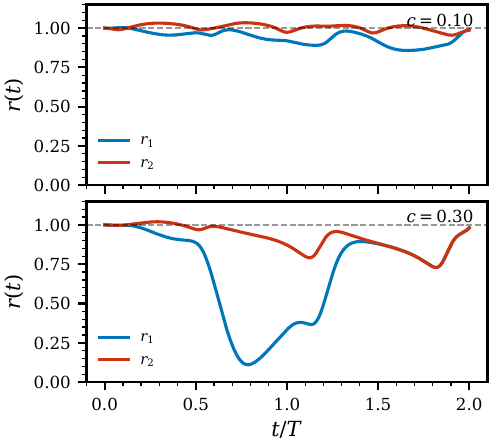}
\caption{Amplitude time histories for on-cycle ($c = 0.10$, top)
and collapse ($c = 0.30$, bottom) regimes. At weak coupling both
oscillators remain near $r = 1$ throughout. At strong coupling
the controller drives oscillator~1 (blue) toward the origin
while oscillator~2 (red) stays on the limit cycle; the natural
restoring dynamics rebuild $r_1$ to above $0.97$ before $T_f$.}
\label{fig:amplitude}
\end{figure}

The picture changes at stronger coupling. At
$c = 0.30$, the controller drives oscillator~1 to
$r_{1,\min} = 0.11$, nearly to the origin, while oscillator~2
remains near the limit cycle. This is clearly visible in the
phase portrait: the trajectory of oscillator~1 spirals inward to
a small neighborhood of the origin, dwells there while
accumulating a phase advantage, then spirals back out to rejoin
the limit cycle in approximate anti-phase with oscillator~2. We
call this \emph{amplitude collapse}.
Fig.~\ref{fig:amplitude} shows the amplitude time histories
$r_1(t)$ and $r_2(t)$ for both regimes: at $c = 0.10$, both
amplitudes remain near unity throughout; at $c = 0.30$,
oscillator~1 collapses to $r_{1,\min} = 0.11$ while oscillator~2
stays on the limit cycle, before the natural dynamics rebuild
$r_1$ to $0.98$ by $T_f$.

The contrast between $c = 0.20$ and $c = 0.30$ is sharp. At $c = 0.20$, the trajectory stays near the limit
cycle ($r_{1,\min} = 0.77$); at $c = 0.30$, it collapses to
$r_{1,\min} = 0.11$. 
Within the coupling values resolved by this sweep, we observe no intermediate regime of moderate amplitude excursion: the computed HJB trajectory changes rapidly from a near-cycle strategy to a deep amplitude excursion.
The dense sweep of Sec.~\ref{sec:transition} resolves the
crossover to $\Delta c < 0.006$ and returns no trajectory with
$r_{1,\min}$ between $0.13$ and $0.72$.
This sharp difference between adjacent
coupling values motivates a dense sweep to precisely locate the
transition.

\subsection{Sharp Strategy Transition}\label{sec:transition}

\begin{figure}[t]
\centering
\includegraphics[width=\columnwidth]{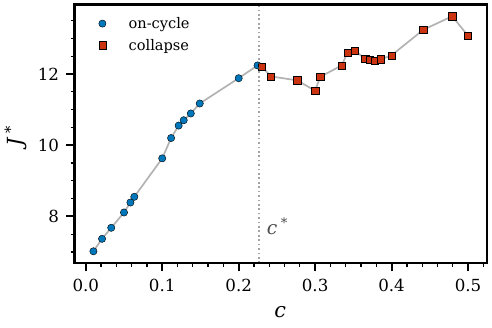}
\caption{Optimal cost $J^*$ versus coupling strength $c$ over 30
values. The strategy transition changes sharply near $c^* \approx 0.227$
($\Delta c < 0.006$). The cost is non-monotonic: it
\emph{decreases} after the transition due to the energetic
advantage of amplitude collapse at this nonisochronicity
($b = 0.8$).}
\label{fig:cost_sweep}
\end{figure}

A dense sweep over 30 coupling values in the range
$c \in [0.01, 0.50]$ reveals a change in the computed HJB strategy at
$c^* = 0.227 \pm 0.003$ (Fig.~\ref{fig:cost_sweep}). At
$c = 0.2247$, the optimal trajectory stays near the limit cycle
with $r_{1,\min} = 0.728$. At $c = 0.2299$, a change of only
$\Delta c = 0.0052$, the optimal trajectory collapses to
$r_{1,\min} = 0.097$. Thus, the change from near-cycle steering to amplitude collapse occurs over an interval narrower than $0.006$.  We therefore use $c^* \approx 0.227$ as an operational estimate of the strategy-change threshold at this resolution. Time at $r_1 < 0.5$ discriminates more sharply than $r_{1,\min}$: zero at $c = 0.2247$, $0.51\,T$ just above it. Repeating the two bracketing solves at $G = 80$ and $G = 96$ places $c^*$ in the same interval, with $r_{1,\min}$ just above threshold varying by $0.007$ and the hold time by $0.03\,T$ across the three grids.
It appears that the controller does not
gradually increase its off-cycle excursion as coupling
increases; instead, it abruptly switches from one strategy to the other.
The optimal cost $J^*(c)$ is also \emph{non-monotonic}:
it decreases between $c = 0.20$ and $c = 0.30$ despite the
stronger coupling opposition (Fig.~\ref{fig:cost_sweep}).
This is a signature of the collapse strategy's
efficiency, which we explain in Sec.~\ref{sec:mechanism}.

\subsection{Physical Mechanism: Phase Repositioning and Phase-Velocity
Surplus}\label{sec:mechanism}

The amplitude-collapse strategy provides two distinct
benefits. First, near the origin the oscillator's phase becomes highly
sensitive to state perturbations, enabling large phase repositioning
without continuously opposing the synchronizing phase dynamics on the
limit cycle. This mechanism remains available at $b=0$. Second, for
$b\neq0$, collapse provides an additional \emph{phase-velocity
surplus}, making the strategy substantially less expensive as
nonisochronicity increases.

For the intrinsic dynamics, $\dot\phi = \omega + br^2\sin 2\phi$.
At $r \approx 0$, the nonisochronous distortion
($\propto b r^2$) vanishes, and the intrinsic angular velocity approaches $\omega$. On the limit cycle ($r = 1$),
the time-averaged angular velocity is
$\Omega = \sqrt{\omega^2 - b^2}$, which is \emph{slower}
than $\omega$ whenever $b \neq 0$. The difference
\begin{equation}\label{eq:delta_omega}
\delta\omega = \omega - \Omega
= \frac{b^2}{\omega + \Omega}
\end{equation}
is the phase-velocity surplus. At $b = 0.8$, $\omega = 1$ gives
$\delta\omega = 0.4$, providing substantial free phase
repositioning.

During the hold interval $\tau$ (the time oscillator~1 spends at
small amplitude), the controller must counter the coupling force
$c\cos\phi_2$ that pushes oscillator~1 away from the origin, at
cost $E_{\mathrm{hold}} \approx Rc^2\tau/2$ (since
$\langle\cos^2\phi_2\rangle = 1/2$), while gaining
$\delta\omega\cdot\tau$ radians of free phase. The trade-off
between this hold cost and the phase gained helps determine whether collapse is energetically favorable.

At $b = 0.8$, $c = 0.30$ the hold costs
$E_{\mathrm{hold}} \approx 2.7$ and returns
$\delta\omega\cdot\tau = 2.4$ radians, so collapse is favorable
when repositioning the same $2.4$ radians on the cycle would cost
more than $2.7$. Amplitude recovery is governed by the radial
dynamics $\dot r = ar(1-r^2)$, which is logistic in $r^2$ rather
than exponential; the Floquet rate $-2a$ describes decay toward the
cycle from nearby and does not apply over the full excursion. Across the
sweeps reported here, amplitude is fully rebuilt by the end of the
horizon, with $r_i(T_f) > 0.97$ in every case.

\subsection{Role of Nonisochronicity}\label{sec:cb_sweep}

Having established the collapse mechanism at $b = 0.8$, we now
investigate how it depends on the nonisochronicity parameter $b$.
Table~\ref{tab:cb_sweep} shows $r_{1,\min}$ across 36
$(c,b)$ pairs.

\begin{table}[t]
\centering
\caption{Minimum amplitude $r_{1,\min}$ across the $(c,b)$
parameter plane. Bold entries indicate amplitude collapse
($r_{1,\min} < 0.3$). All runs use $G = 64$,
$a = \omega = 1$.}
\label{tab:cb_sweep}
\renewcommand{\arraystretch}{1.15}
\begin{tabular}{l cccccc}
\toprule
$b$ & $c{=}0.05$ & $c{=}0.10$ & $c{=}0.20$ & $c{=}0.30$
    & $c{=}0.40$ & $c{=}0.50$ \\
\midrule
0.0 & 0.928 & 0.885 & 0.806
    & \textbf{0.181} & \textbf{0.218} & \textbf{0.234} \\
0.2 & 0.908 & 0.857 & \textbf{0.112}
    & \textbf{0.191} & \textbf{0.219} & \textbf{0.223} \\
0.4 & 0.900 & 0.858 & \textbf{0.113}
    & \textbf{0.190} & \textbf{0.207} & \textbf{0.196} \\
0.6 & 0.901 & 0.862 & \textbf{0.102}
    & \textbf{0.177} & \textbf{0.165} & \textbf{0.272} \\
0.8 & 0.897 & 0.857 & 0.773
    & \textbf{0.110} & \textbf{0.131} & \textbf{0.184} \\
0.9 & 0.870 & 0.845 & 0.801
    & 0.761 & 0.337 & 0.374 \\
\bottomrule
\end{tabular}
\end{table}

First, amplitude collapse occurs at all
$b$ values for sufficiently strong coupling, including the
isochronous case $b = 0$. The estimated strategy-change threshold $c^*$ is
non-monotonic in $b$: it is lowest at
$b = 0.2$ to $0.6$ ($c^* \in (0.10, 0.20)$), higher at $b = 0$
($c^* \in (0.20, 0.30)$), and highest at $b = 0.9$
($c^* > 0.50$). At $b = 0.9$ and $c = 0.40$ the trajectory dips
to $r_{1,\min} = 0.34$ but spends under $0.07\,T$ below $r_1 = 0.5$ out
of a $2\,T$ horizon, so it grazes rather than holds and the
mechanism above does not engage.

Second, while collapse occurs even at $b = 0$, the
additional energetic advantage associated with off-cycle phase accumulation depends critically on
nonisochronicity. At $c = 0.30$, the optimal cost drops from
$J^* = 19.93$ at $b = 0$ to $J^* = 11.53$ at $b = 0.8$, a
$42\%$ reduction. This is because the phase velocity surplus
$\delta\omega$ grows with $b$: at $b = 0$,
$\delta\omega = 0$ and collapse provides no free phase; at
$b = 0.8$, $\delta\omega = 0.4$ and the collapsed
oscillator gains $138^\circ$ of free phase repositioning
(Table~\ref{tab:phase_budget}).

Third, at $b = 0.9$ the long period ($T = 14.4$) provides
sufficient time for on-cycle steering, delaying the onset of
collapse to higher coupling values.

\subsection{Quantitative Verification: Free Phase Budget}

Table~\ref{tab:phase_budget} quantifies the phase accumulation
during the collapse interval (defined as the time oscillator~1
spends at $r_1 < 0.5$) at $c = 0.30$ across nonisochronicity
values. The free phase $\delta\omega\cdot\tau$ matches the
analytical prediction $(\omega - \Omega)\tau$ to within
$1^\circ$ at all $b$ values, confirming the phase velocity
surplus mechanism.

\begin{table}[t]
\centering
\caption{Phase budget during collapse ($r_1 < 0.5$) at
$c = 0.30$. The free phase $\delta\omega\cdot\tau$ is the
phase gained by operating at collapsed amplitude rather than on
the limit cycle.}
\label{tab:phase_budget}
\renewcommand{\arraystretch}{1.15}
\begin{tabular}{c cccc c}
\toprule
$b$ & $\tau/T$ & $\delta\omega\cdot\tau$
    & $\omega\cdot\tau$ & $\Omega\cdot\tau$ & $J^*$ \\
\midrule
0.0 & 0.74 & $0^\circ$    & $268^\circ$ & $268^\circ$ & 19.93 \\
0.2 & 0.71 & $5^\circ$    & $262^\circ$ & $257^\circ$ & 18.10 \\
0.4 & 0.68 & $22^\circ$   & $266^\circ$ & $244^\circ$ & 15.48 \\
0.6 & 0.64 & $58^\circ$   & $288^\circ$ & $230^\circ$ & 14.16 \\
0.8 & 0.57 & $138^\circ$  & $344^\circ$ & $206^\circ$ & 11.53 \\
0.9 & \multicolumn{4}{c}{(no collapse, $r_{1,\min}=0.76$)}
    & 19.76 \\
\bottomrule
\end{tabular}
\end{table}

The optimizer extends the hold time from $0.74\,T$ ($b = 0$) to
$0.57\,T$ ($b = 0.8$) as the per-unit-time phase return increases.
The $42\%$ cost reduction from $b = 0$ to $b = 0.8$ is accounted
for by $138^\circ$ of free phase repositioning and the
consequent reduction in residual steering cost that must be paid
on the limit cycle after amplitude rebuilds.

\section{Comparison with Direct Optimization}\label{sec:direct}

To examine whether direct-optimization methods readily recover the collapse strategy,
we applied both
gradient-based and global stochastic optimization to the same
problem at $c = 0.30$ (Table~\ref{tab:local}). We used
L-BFGS-B~\cite{byrd1995}, a gradient-based quasi-Newton method
with bound constraints, and CMA-ES~\cite{hansen2016}, a
derivative-free evolutionary strategy that adapts a covariance
matrix over a population of candidate solutions. The control
signal was parameterized as piecewise constant over 100 uniform
time intervals.

\begin{table}[t]
\centering
\caption{Direct optimization vs.\ HJB at $c = 0.30$.}
\label{tab:local}
\renewcommand{\arraystretch}{1.15}
\begin{tabular}{lcccc}
\toprule
Method & $J^*$ & $r_{1,\min}$ & Evals & Strategy \\
\midrule
HJB ($G{=}64$) & 11.53 & 0.110 & \textemdash & collapse \\
\midrule
L-BFGS-B (zero init) & 21.23 & 0.746 & \textemdash & on-cycle \\
L-BFGS-B (warm start) & 15.75 & 0.715 & \textemdash & on-cycle \\
L-BFGS-B (5 random) & 27.50 & 0.695 & \textemdash & on-cycle \\
CMA-ES (400D) & 24.24 & 0.964 & 25k & on-cycle \\
CMA-ES (100D) & 24.75 & 0.991 & 196k & on-cycle \\
\bottomrule
\end{tabular}
\end{table}

Neither method discovers the collapse strategy in the runs considered. All L-BFGS-B
runs, from zero initialization, a warm start from an adjacent
coupling value, and five random initializations, converge to
on-cycle solutions with $r_{1,\min} > 0.69$, costing at least
$36\%$ more than HJB. CMA-ES was run with $196{,}000$ function
evaluations and converged to $r_{1,\min} = 0.99$ at $2.1\times$
the HJB cost. 
While these calculations do not rigorously establish that the collapse trajectory belongs to a mathematically distinct basin of attraction, they show that the substantially lower-cost HJB trajectory is not readily recovered by these standard direct-search formulations.

\section{Comparison with Phase-Reduced Control}\label{sec:prc}

Phase reduction is the standard approach for oscillator
control~\cite{moehlis2006,monga2019}. To quantify the
limitations of reduced-order models, we compare the HJB solution
against two phase-reduced baselines.

The \emph{uncoupled baseline}~\cite{moehlis2006} designs
controls independently for each oscillator using the
infinitesimal phase response curve (iPRC), obtained from the
adjoint equation~\cite{brown2004,monga2019}. The \emph{coupled
baseline} incorporates the interaction function
$H(\psi) = (c/T)\int_0^T Z_x(\phi)\,[x_\gamma(\phi{+}\psi) -
x_\gamma(\phi)]\,d\phi$
from~\cite{ermentrout2010} and solves the resulting
two-dimensional optimal control problem via Pontryagin shooting.
Both baselines are applied open-loop to the full four-dimensional
system.

\begin{table}[t]
\centering
\caption{Three-way comparison: uncoupled PRC, coupled PRC,
and full-state HJB.}
\label{tab:prc_vs_hjb}
\renewcommand{\arraystretch}{1.15}
\begin{tabular}{c ccc ccc}
\toprule
 & \multicolumn{3}{c}{$\Delta\phi(T_f)$}
 & \multicolumn{3}{c}{$J_u$} \\
\cmidrule(lr){2-4}\cmidrule(lr){5-7}
$c$ & Uncoup. & Coupled & HJB
    & Uncoup. & Coupled & HJB \\
\midrule
0.01 & $178.3^\circ$ & $175.2^\circ$ & $179.7^\circ$
     & 7.47 &  8.14 &  7.02 \\
0.05 & $170.8^\circ$ & $177.6^\circ$ & $179.9^\circ$
     & 7.47 &  9.28 &  8.11 \\
0.10 & $140.3^\circ$ & $178.6^\circ$ & $180.0^\circ$
     & 7.47 & 10.95 &  9.63 \\
0.20 &  $48.7^\circ$ & $171.5^\circ$ & $179.9^\circ$
     & 7.47 & 15.07 & 11.88 \\
0.30 &  $23.7^\circ$ & $171.4^\circ$ & $180.0^\circ$
     & 7.47 & 20.01 & 11.53 \\
0.50 &  $10.9^\circ$ & $150.9^\circ$ & $180.0^\circ$
     & 7.47 & 30.54 & 13.07 \\
\bottomrule
\end{tabular}
\end{table}

\begin{figure}[t]
\centering
\includegraphics[width=\columnwidth]{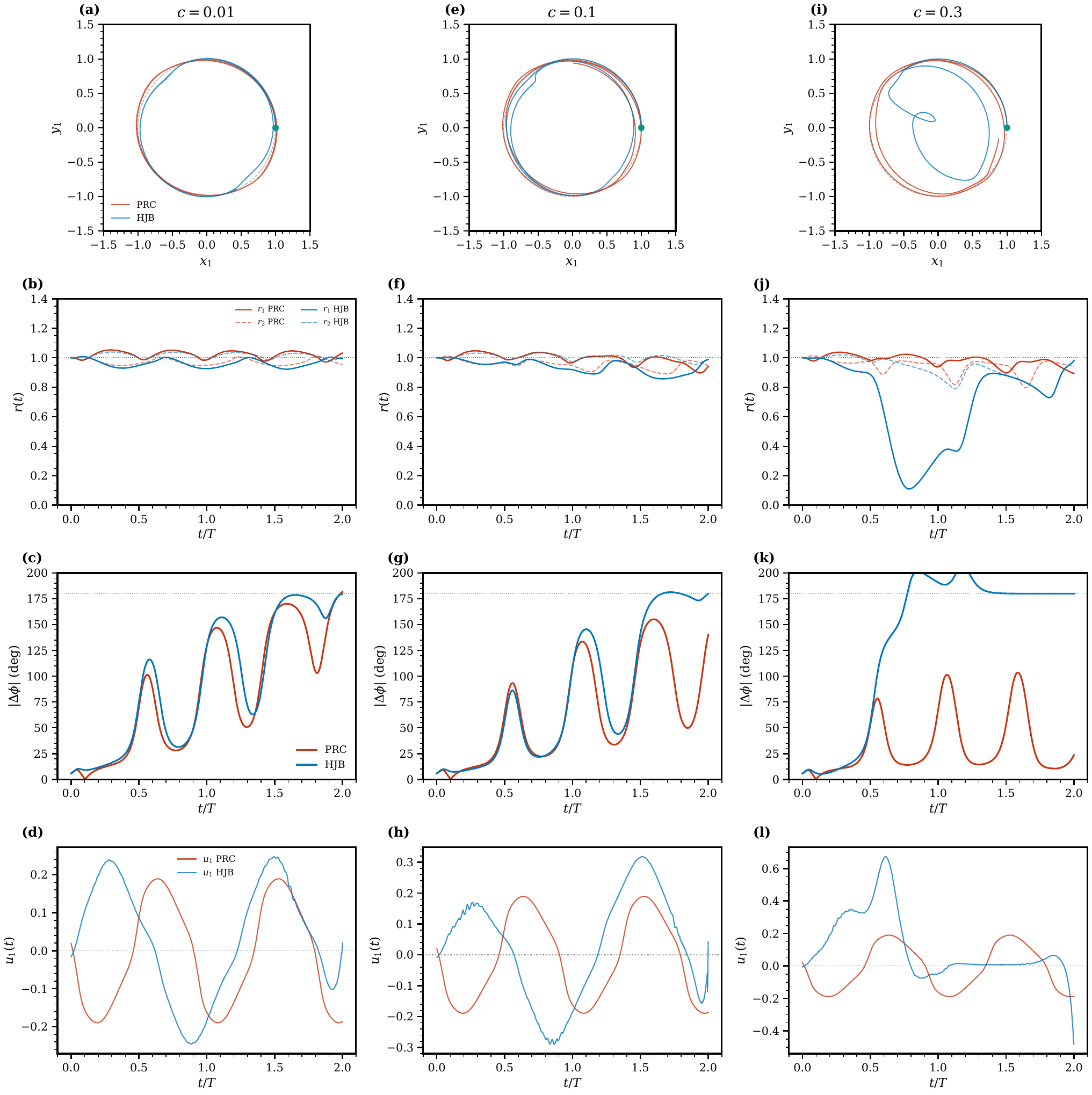}
\caption{PRC (red) vs.\ HJB (blue) optimal trajectories for
$c = 0.01$ (left), $c = 0.10$ (center), and $c = 0.30$ (right).
At weak coupling, both methods produce similar near-limit-cycle
trajectories. At $c = 0.30$, the HJB controller collapses
amplitude while PRC remains constrained to the limit cycle.}
\label{fig:prc_vs_hjb}
\end{figure}

Table~\ref{tab:prc_vs_hjb} and Fig.~\ref{fig:prc_vs_hjb}
summarize the comparison. At weak coupling ($c = 0.01$), all
three methods perform comparably: the phase difference reaches
$175$ to $180^\circ$ at similar cost, and the HJB trajectory
stays near the limit cycle. The coupled phase model improves
significantly over the uncoupled baseline at moderate coupling:
at $c = 0.10$, coupled PRC achieves $178.6^\circ$ versus
$140.3^\circ$ uncoupled. However, the coupled model saturates
near $171^\circ$ for $c \ge 0.20$ and cannot reach the target,
while costing $1.3$ to $2.3$ times more energy than HJB.

The discrepancy at strong coupling is not only quantitative. The HJB trajectory at $c=0.30$ reaches $r_{1,{\rm min}} = 0.11$, whereas the phase-only models constrain the state to the limit-cycle manifold and therefore cannot represent this mechanism. Phase-amplitude reductions can incorporate transverse dynamics and extend beyond standard weak-input phase reduction~\cite{wilson2020adapt,dewanjee2024,wilson2019,monga2020,toth2025}. However, the collapse trajectory considered here, with $|r - 1| =0.89$, lies considerably deeper off-cycle than the regimes in which the cited reduced-order control methods have been demonstrated.


\section{Validation on Van der Pol Oscillators}
\label{sec:vdp}

To test whether off-cycle optimal control is specific to the
Stuart-Landau-like algebraic structure, we solved the 4D HJB
equation for two coupled Van~der~Pol oscillators:
\begin{equation}
\dot x_i = y_i + c(x_j{-}x_i) + u_i, \qquad
\dot y_i = \mu(1{-}x_i^2)y_i - x_i,
\end{equation}
with $\mu = 1$ and $\umax = 0.5$. The Van~der~Pol limit cycle
is non-circular (egg-shaped) with amplitude $\approx 2$ in $x$
and period $T = 6.66$, distinguishing it qualitatively from the
unit-circle Stuart-Landau normal form. 
Unlike the latter system, the Van~der~Pol oscillator does not admit the same simple decomposition of the angular velocity into an on-cycle mean frequency and an analytically defined off-cycle phase
velocity surplus. We therefore use this example primarily to test whether strong coupling can induce the HJB controller to exploit large departures from the limit cycle in a qualitatively different oscillator.

\begin{table}[t]
\centering
\caption{Van~der~Pol coupling sweep ($\mu = 1$, $G = 64$).
Off-cycle deviation is measured as
$\max_t\min_s\|z_i(t) - x_\gamma(s)\|$.}
\label{tab:vdp}
\renewcommand{\arraystretch}{1.15}
\begin{tabular}{c c cc c}
\toprule
$c$ & $J^*$ & Off-cycle 1 & Off-cycle 2 & $\Delta\phi$ \\
\midrule
0.05 & 25.27 & 0.52 & 0.51 & $180.0^\circ$ \\
0.10 & 12.09 & 0.39 & 0.39 & $179.9^\circ$ \\
0.15 & 18.77 & 0.41 & 0.58 & $179.9^\circ$ \\
0.20 & 27.96 & 0.62 & 0.80 & $179.5^\circ$ \\
0.25 & 37.76 & 1.09 & 0.86 & $180.0^\circ$ \\
0.30 & 48.18 & 1.20 & 0.91 & $179.5^\circ$ \\
\bottomrule
\end{tabular}
\end{table}

For weak to moderate coupling ($c \le 0.20$), both oscillators
remain close to the limit cycle, with off-cycle deviations below
$0.80$ (less than $40\%$ of the limit-cycle amplitude). At
$c = 0.25$, the optimal trajectory deviates by $1.09$ from the
limit cycle ($55\%$ of the limit-cycle amplitude), increasing to
$1.20$ at $c = 0.30$. These significant off-cycle excursions
confirm that the HJB controller exploits the full state space
even for oscillators with qualitatively different limit-cycle
geometry. 
Thus, large off-cycle excursions selected by full-state optimal control are not specific to the Stuart-Landau-like algebraic structure.

\section{Discussion}\label{sec:discussion}

The amplitude-collapse mechanism reflects two distinct geometric effects. First, near the origin the oscillator's phase is highly sensitive to state perturbations, enabling efficient phase repositioning without continuously opposing the synchronizing phase dynamics on the limit cycle. This effect can favor collapse even at $b=0$.  Second, for $b\neq0$, the phase-velocity surplus $\delta\omega=\omega-\Omega>0$ provides an additional energetic benefit by allowing faster intrinsic phase accumulation at small amplitude.

The closest precedent for using HJB to exploit off-cycle geometry
in oscillator control is~\cite{danzl2009}, which drives a single
oscillator to its phaseless set for phase randomization. The
present work differs in both objective and setting: we target a
specific phase relation (anti-phase) in a coupled pair, and the
off-cycle mechanism, transient amplitude collapse to neutralize
coupling and harvest phase drift, serves a constructive
synchronization purpose rather than a disruptive randomization
one.

Several papers in the phase-amplitude
literature~\cite{wilson2020adapt,wilson2020b,dewanjee2024}
explicitly aim to extend validity beyond the weak-perturbation
regime. The present result does not claim that such extensions are
impossible in principle for deep excursions. Rather, we observe
that within the reviewed literature, no reduced-order control
study discovers or certifies an optimal strategy involving
excursions as deep as $|r-1| = 0.89$, and the collapse trajectory
lies well outside the regime in which these methods have been
validated. The full-state HJB calculation demonstrates that a substantially lower-cost control mechanism can involve excursions far outside the limit-cycle neighborhood on which standard phase reduction is based. 

Because the anti-phase orbit is unstable at these couplings,
the achieved configuration does not persist once control is removed:
for $c=0.2299$, the phase difference stays above $170^\circ$ for
roughly two periods and decays to in-phase within about ten.
Maintenance is nonetheless cheap. The anti-phase subspace
$\{x_2=-x_1,\,y_2=-y_1\}$ is invariant, requiring no control without
disturbances; under additive state noise of standard deviation $0.01$,
feedback damping the in-phase component maintains anti-phase for about
$2\%$ of $J^*$ per horizon. Thus, desynchronization may be expensive
to establish but inexpensive to sustain under continued actuation.

The practical implication is not that a physical
oscillator should be driven near zero amplitude, but that temporary
off-cycle excursions may reduce the cost of overcoming strong
synchronizing coupling. This motivates full-state or phase-amplitude
control when near-cycle steering becomes energetically expensive.

Extension to systems of three or more oscillators ($\geq 6$D) is
limited by the curse of dimensionality inherent in grid-based HJB
solvers. Tensor, neural, and actor--critic approximations offer
alternatives to full grids~\cite{dolgov2021,oster2022,tassa2007,nakamura2019,zhou2021};
whether they preserve the sharp strategy transition in
Sec.~\ref{sec:transition} merits study. Related surrogate and
model-reduction methods learn reduced representations of expensive
physical models~\cite{gao2020,fani2026,mohammadagha2025,mazloom2026},
while synthetic data, learned priors, and model compression provide
complementary tools~\cite{fathijam2026,saberi2026,madinei2026}.
These ideas motivate reduced-cost approximations of offline HJB
feedback maps. Optimization under physical constraints also appears
in production, mapping, sensing, vibration, and actuator-control
problems~\cite{sarrafan2023,orisatoki2026,moradi2026,hashemi2025,jalalvand2016}.

\section{Conclusion}\label{sec:conclusion}

We solved the four-dimensional HJB equation for two coupled Stuart-Landau-like oscillators and found a sharp change in the numerically optimal control strategy as coupling increases. At weak coupling, the HJB trajectory remains near the limit cycle and is well described by phase-based intuition. At stronger coupling, it instead makes a deep off-cycle excursion in which one oscillator's amplitude is transiently reduced before returning to the target anti-phase orbit. The results identify two benefits of this strategy: reducing the influence of synchronizing coupling, which can favor collapse even for isochronous dynamics, and exploiting a phase-velocity surplus at small amplitude, which provides an additional energetic advantage in the nonisochronous case. Numerical comparisons show that the particular phase-reduced baselines considered here do not reproduce this deep off-cycle mechanism, while the tested gradient-based and stochastic direct optimizers do not recover the lower-cost HJB trajectory. Results for coupled Van der Pol oscillators further indicate that exploitation of off-cycle dynamics for optimal control is not specific to the Stuart-Landau-like model. More broadly, these computations suggest that when coupling is sufficiently strong, optimal oscillator control may benefit from exploiting regions of state space far from the nominal limit cycle.



\end{document}